\documentclass[12pt]{article}

\usepackage{amsmath}
\usepackage{amsthm}
\usepackage{amsfonts}
\usepackage{amssymb}
\usepackage{graphicx}

\theoremstyle{plain}
\newtheorem{theorem}{Theorem}
\newtheorem{lemma}{Lemma}
\newtheorem{proposition}{Proposition}
\newtheorem{definition}{Definition}
\newtheorem{corollary}{Corollary}

\theoremstyle{definition}
\newtheorem{remark}{Remark}

\title{Discrete L\"owner evolution}
\author{Robert O. Bauer\\ Department of Mathematics\\ University of Illinois at Urbana-Champaign\\ 1409 West Green Street \\ Urbana, IL 61801, USA\\
rbauer@math.uiuc.edu}

\begin{document}

\maketitle

\begin{abstract}
	We study a one parameter family of discrete L\"owner 	evolutions driven by a random 	walk on the real line. We show 	that it converges to the 	stochastic L\"owner evolution (SLE) 	under rescaling. We show that the discrete L\"owner evolution 	satisfies Markovian-type and symmetry properties analogous 	to SLE, and establish a phase transition property for the 	discrete L\"owner evolution when the parameter equals 4.   
\end{abstract}

\section{Introduction}

In this paper we study a discrete version of the stochastic L\"owner evolution ($\text{SLE}_{\kappa}$) introduced by O. Schramm in \cite{schramm:2000}. Whereas SLE is driven by a one dimensional Brownian motion, our discrete L\"owner evolution is driven by a random walk. SLE is a one parameter family of processes of growing random sets in a domain in the plane. We will only consider chordal SLE and our discrete version, where the random sets grow in the upper half-plane from 0 to $\infty$.

It has been shown that, in a sense that can be made precise (\cite{lawler.werner:2000}), any random process of growing sets in the plane that satisfies a certain Markovian type property is given by $\text{SLE}_{\kappa}$ for some $\kappa\in[0,\infty)$. Since SLE is amenable to computations this  led to some spectacular calculations of various quantities long believed out of reach for mathematicians. For example, in a sequence of papers \cite{lsw:2001i},\cite{lsw:2001ii},\cite{lsw:2002},\cite{lsw:2000}, Lawler, Schramm, and Werner calculated all intersection exponents for Brownian motion in the plane. Many of these exponents had been predicted by physicists based on non-rigorous methods from conformal field theory. In particular, Lawler, Schramm, and Werner confirmed a conjecture of Mandelbrot, that the Brownian frontier has Hausdorff dimension 4/3. Furthermore, SLE has been shown to be the scaling limit of various discrete systems, e.g. loop-erased random walk and the outer boundary of critical percolation clusters on the triangular lattice, and is conjectured to give the scaling limit of others, such as the self-avoiding random walk. To confirm the conjectures the existence of the scaling limit and the conformal invariance of the scaling limit need to be established, the latter usually being the main obstacle. 

In this paper we study a discrete (in time) approximation of SLE. Instead of a continuous family of conformal maps $\{f_t\}_{t\in[0,\infty)}$ from the upper half-plane $\mathbb H$ into $\mathbb H$ so that $f_t(\mathbb H)\supseteq f_s(\mathbb H)$ if $t\le s$, we consider a sequence $\{f(m)\}_{m=0}^{\infty}$ of such maps. The ``increments'' $f_m^{-1}\circ f_{m+1}$ are all of the form
\[
	z\in\mathbb H\mapsto S(m)+\sqrt{(z-S(m))^2-4}\in\mathbb H,
\]
where $\{S(m)\}_{m=0}^{\infty}$ is a random walk on $\mathbb R$ with centered increments of variance~$\kappa$. We show in Theorem \ref{T:dle} that the law of  $\{f(m)\}_{m=0}^{\infty}$, properly rescaled, converges weakly to $\text{SLE}_{\kappa}$. The proof relies on Donsker's invariance principle and continuity properties of L\"owner's differential equation, considered as a map from piecewise continuous curves $\psi$ to 1-parameter families of conformal maps $f:\mathbb H\to\mathbb H$. 

To establish continuity, we first choose a topology on the space of conformal maps $f:\mathbb H\to\mathbb H$. One natural choice is the topology of uniform convergence on compacts. In fact, in our context, this is equivalent to uniform convergence on $\{z\in\mathbb H:\Im(z)>a\}$ for every $a\in(0,\infty)$. However, the regularity of the L\"owner equation allows us to choose a stronger topology that also takes the boundary behavior of $f$ into account. We introduce this topology in the context of Cauchy transforms of probability measures in Lemmas  \ref{L:Polish}, \ref{L:reciprocal}, and \ref{L:metric}.    

In Section \ref{S:dle} we study properties of the the sequence $\{f(m)\}_{m=0}^{\infty}$. We show in Theorem~\ref{T:iid} that if the increments $S(m+1)-S(m)$ have the appropriate properties, then $\{f(m)\}_{m=0}^{\infty}$ has the same Markovian-type and symmetry properties as SLE. We call $\{f(m)\}_{m=0}^{\infty}$ a discrete L\"owner evolution with parameter $\kappa$ ($\text{DLE}_{\kappa}$), if the increments $S(m+1)-S(m)$ are centered, independent and identically distributed random variables with variance $\kappa$.  Next, we study the dependency of the discrete L\"owner evolution on $\kappa$. In the paragraphs following Proposition~\ref{P:connect} we describe, graphically, DLE in the special case when the increments $S(m+1)-S(m)$ are Bernoulli random variables. The behavior of the omitted set, i.e. the image of $\mathbb H$ under $f(m)$, is rather easily understood in terms of the underlying random walk $\{S(m)\}$. In our view this connection is not as apparent in the continuous case and making it more explicit is one of our motivations for this paper.  In Proposition~\ref{P:connect} we note the transition from connected to disconnected complement of the image at $\kappa=4$. 
In Theorem~\ref{T:Bessel} we show that Markov chains (with uncountable state space) naturally associated to DLE have a transition from transient to recurrent at $\kappa=4$. These  Markov chains are the discrete analogues of Bessel processes naturally occurring in the study of SLE \cite{werner:2002}.

Finally, we collect in the appendix some facts about monotonic independence in noncommutative probability and its relation to the (deterministic) L\"owner evolution. The impetus to build a discrete L\"owner evolution from the maps $z\mapsto a+\sqrt{(z-a)^2-4}$ came from the preprint \cite{muraki:2000}, which H. Bercovici had kindly brought to our attention.

\section{A discrete approximation of SLE}

Denote ${\mathfrak P}(\mathbb R)$ the space $C([0,\infty);\mathbb R)$ of continuous paths $\psi:[0,\infty)\to\mathbb R$ and endow ${\mathfrak P}(\mathbb R)$  with the topology of uniform convergence on compact intervals. Let $\{X_n\}_{n=1}^{\infty}$ be a sequence of independent real-valued random variables on a probability space $(\Omega,{\cal F},P)$, and assume that the $X_n$'s have mean-value 0, variance $\kappa>0$ and satisfy
\[
	\lim_{R\to\infty}\sup_{n\in\mathbb Z^+}
	\mathbb E^P[|X_n|^2,|X_n|\ge R]=0.
\]
Next, for $n\in\mathbb Z^+$, define $\omega\in\Omega\mapsto S_n(\cdot,\omega)\in{\mathfrak P}(\mathbb R)$ so that $S_n(0,\omega)=0$ and, for each $m\in\mathbb Z^+$, $S_n(\cdot,\omega)$ is linear on the interval $[\frac{m-1}{n},\frac{m}{n}]$ with slope $n^{1/2} X_m(\omega)$. That is,
\[
	S_n(0,\omega)=0,\quad
	S_n\left(\frac{m}{n},\omega\right) 
	=n^{-1/2}\sum_{k=1}^m X_k,\ m\in\mathbb Z^+
\]
and 
\[
	S_n(t,\omega)=(m-nt) S_n\left(\frac{m-1}{n},\omega\right)
	+(1-(m-nt))S_n\left(\frac{m}{n},\omega\right)
\]
for $t\in(\frac{m-1}{n},\frac{m}{n})$. Finally, let 
\[
	\mu_n\equiv(S_n)_* P
\]
denote the distribution of 
$\omega\in\Omega\mapsto S_n(\cdot,\omega)\in{\mathfrak P}(\mathbb R)
$
under $P$. Then it is well known, see \cite{stroock:1993}, that $\mu_n\Longrightarrow\cal W_{\kappa}$ as $n\to\infty$, where $\cal W_{\kappa}$ is the distribution of $\psi\in\mathfrak P(\mathbb R)\mapsto \sqrt{\kappa}\psi\in\mathfrak P(\mathbb R)$ under Wiener's measure $\cal W$ on $\mathfrak P(\mathbb R)$.  

Denote $\mathbb H$ the upper half-plane $\{z\in\mathbb C:\Im(z)>0\}$. 
For $z\in\mathbb H$,  and $\psi\in{\mathfrak P}(\mathbb R)$ consider the chordal L\"owner equation
\begin{equation}\label{E:loewner}
	\frac{\partial}{\partial t}g(t,\psi;z) 	=\frac{2}{g(t,\psi;z)-\psi(t)},\quad 	g(0,\psi;z)=z.
\end{equation}
Then
\begin{equation}\label{E:exbound}
	\left|\frac{\partial}{\partial t}g(t,\psi;z)\right|\le\frac{2}{{\Im}(
	g(t,\psi;z))},
\end{equation}
and
\begin{equation}\label{E:ibound}
	\frac{\partial }{\partial t}{\Im}(g(t,\psi;z))
	=-\frac{2\Im(g(t,\psi;z))}{|g(t,\psi;z)-\psi(t)|^2}<0.
\end{equation}
The inequalities imply in particular that for each $z\in\mathbb H$ and $\psi\in{\mathfrak P}(\mathbb R)$ the solution is well defined up to a time $\tau(\psi;z)\in(0,\infty]$, and that if $\tau(\psi;z)<\infty$, then $\lim_{t\nearrow\tau(\psi;z)}{\Im}(g(t,\psi;z))=0$. Let ${\cal K}(t,\psi)$ be the closure of $\{z\in\mathbb H:\tau(\psi;z)\le t\}$.

\begin{proposition}[\cite{lawler:2001}]\label{P:lawler} For every $t\in(0,\infty)$ and 	$\psi\in\mathfrak P(\mathbb R)$, $g(t,\psi;\cdot)$ is a conformal 	transformation of $\mathbb H\backslash{\cal K}(t,\psi)$ onto 	$\mathbb H$ satisfying
	\[
		g(t,\psi;z)=z+\frac{2t}{z}+O\left(\frac{1}{|z|^2}\right),\quad z\to\infty.
	\]
\end{proposition}

Let ${\bf M}_1(\mathbb R)$ be the set of Borel probability measures on $\mathbb R$ and denote ${\bf M}_1^0(\mathbb R)$ the subset of Borel probability measures with compact support. If $\mu\in{\bf M}_1^0(\mathbb R)$ let $[A_{\mu},B_{\mu}]$ denote the convex closure of $\text{supp}(\mu)$. 
For $\mu,\nu\in {\bf M}_0^1(\mathbb R)$ let 
\[
	\rho(\mu,\nu)=\rho_L(\mu,\nu)
	+\max\{|A_{\mu}-A_{\nu}|,|B_{\mu}-B_{\nu}|\},
\]
where
\begin{align}
	\rho_L\equiv\inf\{\delta:\mu((-\infty,x-\delta])-\delta&\le\nu	((-\infty,x])\notag\\
	&\le\mu((-\infty,x+\delta])+\delta\text{ for all }x\in\mathbb R\}
	\notag
\end{align}
is the L\'evy distance between $\mu$ and $\nu$.
Then $\rho$ is  a metric on ${\bf M}_1^0(\mathbb R)$. 

\begin{lemma}\label{L:Polish} $({\bf M}_1^0(\mathbb R),\rho)$ is a 	Polish space
\end{lemma}

\begin{proof}
If $\{\mu_n\}_{n=1}^{\infty}\subset{\bf M}_1^0(\mathbb R)$ is a $\rho$-Cauchy sequence, then $\{\mu_n\}_{n=1}^{\infty}$ is also a $\rho_L$-Cauchy sequence in ${\bf M}_1(\mathbb R)$ and 
$\overline{\bigcup_{n=1}^{\infty}\text{supp}(\mu_n)}$
is compact. Since ${\bf M}_1(\mathbb R)$ is $\rho_L$-complete there is a $\mu\in {\bf M}_1(\mathbb R)$ so that $\rho_L(\mu_n,\mu)\to0$ as $n\to\infty$. Since also $\text{supp}(\mu)\subset\overline{\bigcup_{n=1}^{\infty}\text{supp}(\mu_n)}$, $\mu\in{\bf M}_1^0(\mathbb R)$ and ${\bf M}_1^0(\mathbb R)$ is $\rho$-complete. Finally, it is easy to see that the set of all convex combinations $\sum_{k=1}^n \alpha_k\delta_{x_k}$, where $n\in\mathbb Z^+$, $\{\alpha_k\}_{k=1}^n\subset[0,1]\cap\mathbb Q$ with $\sum_{k=1}^n\alpha_k=1$, and $\{x_k\}_{k=1}^n\subset\mathbb Q$, is a countable $\rho$-dense set in ${\bf M}_1^0(\mathbb R)$.
Thus $({\bf M}_1^0(\mathbb R),\rho)$ is a Polish space.
\end{proof}

Given $\mu\in{\bf M}_1(\mathbb R)$, denote $G_{\mu}$ its Cauchy transform 
\[
	z\in\mathbb C\backslash\text{supp}(\mu)\mapsto 	G_{\mu}(z)=\int_{\mathbb R}
	\frac{\mu(dx)}{z-x}\in\mathbb C\backslash\{0\}.
\]
Note that $G_{\mu}$ is analytic, and that $G_{\mu}(\bar{z})=\overline{G_{\mu}(z)}$. Furthermore, $G_{\mu}$  cannot be extended analytically beyond $\mathbb C\backslash\text{supp}(\mu)$. Indeed, if $G_{\mu}$ extends analytically to $x\in\mathbb R$ then it must extend to a neighborhood $(x-\delta, x+\delta)$ for some $\delta>0$. By continuity we then have $\lim_{y\searrow0}\Im(G(a+\sqrt{-1}y))=0$, uniformly on compact subsets of $(x-\delta, x+\delta)$. But by Stieltjes' inversion formula \cite[2.20]{deift:2000}, for any $-\infty<x<x'<\infty$
\begin{equation}
	\frac{1}{\pi}\lim_{y\searrow0}\int_x^{x'}\Im(G_{\mu}
	(a+\sqrt{-1}y)
	)\ da=\mu((x,x'))+\frac{1}{2}[\mu(\{x\})+\mu(\{x'\})].
\end{equation}
Hence $\mu((x-\delta,x+\delta))=0$ and $x\notin\text{supp}(\mu)$.

Since $G_{\mu}(z)\neq0$ for all $z\in\mathbb C\backslash\text{supp}(\mu)$ we may define the reciprocal Cauchy transform $f_{\mu}:\mathbb H\to\mathbb H$ by $f_{\mu}(z)=1/G_{\mu}(z)$.

\begin{lemma}\label{L:reciprocal}
	An analytic function $f:\mathbb H\to\mathbb H$ is the 	reciprocal  Cauchy transform of some compactly supported 	probability measure $\mu$ on $\mathbb R$, if and only if 
	\begin{equation}\label{E:reciprocal}
		\inf_{z\in\mathbb H}\frac{\Im(f(z))}{\Im(z)}=1
	\end{equation}
	and $G=1/f$ extends analytically to $\mathbb C\backslash	[-N,N]$ for some $N\in\mathbb N$. 	
\end{lemma} 

\begin{proof}
	By \cite[Proposition 2.1]{maassen:1992}, $f$ is the reciprocal 	Cauchy transform of a probability measure $\mu$ on 	$\mathbb R$ if and only if \eqref{E:reciprocal} holds. If $\mu$ 	has compact support $K\subset\mathbb R$ and $f=f_{\mu}$, 	then $G=1/f$ extends to $\mathbb C\backslash K$.  	Conversely, if $f$ satisfies \eqref{E:reciprocal} and $G=1/f$ 	extends analytically to $\mathbb C\backslash[-N,N]$, then by  	\cite[Proposition 	2.1]{maassen:1992}, $f=f_{\mu}=1/G_{\mu}$  	for some $\mu\in{\bf M}_1(\mathbb R)$, and then by Stieltjes' 	inversion formula $\text{supp}(\mu)\subset[-N,N]$. 
\end{proof}

Recall that for a domain $D\subseteq\mathbb C$ a function $f:D\to\mathbb C$ is univalent if it is analytic and 1-1. Let 
\[
	{\bf M}^U=\{\mu\in{\bf M}_1^0(\mathbb R):G_{\mu}:\mathbb 	H\to\mathbb C\backslash\overline{\mathbb H}
	\text{ is univalent}\}
\]
If $f:\mathbb H\to\mathbb H$ is univalent, then we may extend $f$ to $\mathbb C\backslash\mathbb R$ as a univalent function by the Schwarz reflection principle. We say $f$ has a univalent extension to $\mathbb C\backslash[a,b]$, if $f$ extends as a univalent function to $\mathbb C\backslash[a,b]$. Finally, define $A_f,B_f\in\mathbb R$ by
\[
	[A_f,B_f]=\bigcap\{[a,b]:f \text{ has a univalent extension to 
	}\mathbb C\backslash[a,b]\},
\]
whenever the right-hand side is nonempty.

\begin{lemma}\label{L:metric}
	If $\mu\in{\bf M}^U$, $A_{\mu}\neq B_{\mu}$ and $f=f_{\mu}$, 	then $[A_f,B_f]=[A_{\mu},B_{\mu}]$. 
	Furthermore,  $({\bf M}^U,\rho)$ is a Polish space.
	Finally, if 	$\{\mu_n\}_{n=1}^{\infty}\cup\{\mu\}\subset 	{\bf M}^U$ and $f=f_{\mu}$, $f_n=f_{\mu_n}$, $n\in\mathbb 	Z^+$, then $\rho(\mu_n,\mu)\to0$, as 	$n\to\infty$, if and only 	if $f_n\to f$ uniformly on $\{z\in\mathbb C:\Im(z)>a\}$ for any 	$a\in(0,\infty)$, and $\max\{|A_f-A_{f_n}|,|B_f-B_{f_n}|\}$ 	converges to $0$, as $n\to\infty$. 
\end{lemma}
	
\begin{proof}	
	If $\mu\in{\bf M}^U$ and $G=G_{\mu}$, then it is easy to see 	that $[A_G,B_G]=[A_{\mu},B_{\mu}]$. Indeed, since $G_{\mu}$ 	cannot be extended analytically beyond $\mathbb 	C\backslash\text{supp}(\mu)$, it is clear that 	$[A_G,B_G]\supseteq[A_{\mu},B_{\mu}]$. Furthermore, 	$G((B_{\mu},\infty))\subseteq\mathbb R^+$ and 	$G\upharpoonright(B_{\mu},\infty)$ is strictly decreasing. 	Similarly, $G((-\infty,A_{\mu}))\subseteq\mathbb R^-$ and 
	$G\upharpoonright(-\infty,A_{\mu})$ is strictly increasing. 	Hence $[A_G,B_G]=[A_{\mu},B_{\mu}]$. Finally, since $G$ 	does not assume the value zero it follows that $f$ extends as a 	univalent function to $\mathbb C\backslash[A_G,B_G]$. 	Now note that $f$, on a domain of univalence, can only 	assume the value zero once. Since $A_G\neq B_G$ we get 	$[A_f,B_f]=[A_G,B_G]$.
    
	For the following statements, we begin by checking that weak 	convergence of a sequence 	$\{\mu_n\}_{n=1}^{\infty}\subset{\bf M}_1(\mathbb R)$ to a 	probability measure $\mu$ is equivalent to the uniform 	convergence $G_{\mu_n}(z)\to G_{\mu(z)}$ on $\{z\in\mathbb 	C:\Im(z)>a\}$ for any $a>0$. By 
	\cite[Theorem 2.5]{maassen:1992}, 	$\mu_n\Longrightarrow\mu$ as $n\to\infty$, if and only if there 	exists $y>0$ such that 
	\[
		\lim_{n\to\infty}G_{\mu_n}(x+\sqrt{-1}y))
		=G_{\mu}(x+\sqrt{-1}y),\quad x\in\mathbb R.
	\]
	To complete this part assume now that 	$\mu_n\Longrightarrow\mu$ as $n\to\infty$. Then 
	\[\lim_{n\to\infty}G_{\mu_n}(z)=G_{\mu z}, \quad z\in\mathbb 	H.
	\]
	 Now note that $|G_{\nu}(z)|\le1/\Im(z)$ for all $z\in\mathbb 	H$, and all $\nu\in{\bf M}_1(\mathbb R)$. Hence the family 	$\{G_{\nu},\nu\in{\bf M}_1(\mathbb R)\}$ is locally bounded in 	$\mathbb H$ and it follows from Vitali's theorem, 	\cite{duren:1983}, that $G_{\mu_n}(z)\to G_{\mu}(z)$ uniformly 	on compacts. Since $\mu_n\Longrightarrow\mu$ as 	$n\to\infty$, for any $\epsilon>0$ there exists $N>0$ so that 
	$\sup_n \mu_n(\mathbb R\backslash[-N,N])\le\epsilon$.  	Hence $|G_{\mu_n}(z)|\le 1/N+\epsilon/a$ on $\{z\in\mathbb C: 	\Im(z)>a\text{ and }|z^2|>2N\}$ and it follows that 	$G_{\mu_n}(z)\to G_{\mu}(z)$ uniformly on $\{z\in\mathbb 	C:\Im(z)>a\}$. Since uniform (on compacts) limits of univalent 	functions are either univalent or constant (\cite{duren:1983}), 	and since $G_{\mu}$ cannot be constant as $G_{\mu}(z)\to0$ 	as $z\to\infty$, it follows that ${\bf M}^U$ is $\rho$-closed 	and this implies the second statement.  

	Next, given a compact $A\subset\mathbb H$, 
	\[
		d\equiv\inf_{z\in A}|G_{\mu}(z)|>0
	\]
	and there is an $N\in\mathbb Z^+$ such that for all $n\ge 	N$, $\inf_{z\in A}|G_{\mu_n}(z)|\ge d/2$. Hence, for $n\ge N$,
	\begin{align}
		\sup_{z\in A}|f_{\mu}(z)-f_{\mu_n}(z)|&= \sup_{z\in A}
		\frac{|G_{\mu_n}(z)-G_{\mu}(z)|}{|G_{\mu}(z)G_{\mu_n}(z)|}
		\notag\\
		&\le\frac{2}{d^2}\sup_{z\in A}|G_{\mu_n}(z)-G_{\mu}(z)|
		\to0,\notag
	\end{align}
	as $n\to\infty$. Since $\overline{\bigcup_{n=1}^{\infty}
	\text{supp}(\mu_n)}$ is compact it follows in particular that the 	mean values $\{m_n\}_{1=1}^{\infty}$ of 	$\{\mu_n\}_{n=1}^{\infty}$ converge to the mean value	$m$ of $\mu$ and  from Taylor's formula that there exists an 	$N>0$ and a function $c_n(z)$ such that 	$\sup_{|z|>N}\sup_n|c_n(z)|<\infty$ and so that	\[
		G_{\mu_n}(z)=\frac{1}{z}+\frac{m_n}{z^2}
		+\frac{c_n(z)}{z^3}.
	\]
	This implies that $f_n(z)=z-m_n+e_n(z)/z$, where 	$\sup_n|e_n(z)|$ is uniformly bounded for $|z|>N'$ for some 	$N'>0$. Together with the uniform convergence on compacts 	this gives the uniform convergence of $\{f_n\}_{n=1}^{\infty}$ 	on $\{z\in\mathbb C:\Im(z)>a\}$ for any $a>0$. 
\end{proof}

\begin{remark} Based on the above proof it is easy to show that 
	$\rho(\mu_n,\mu)\to0$, as 	$n\to\infty$, if and only if for every 	$\epsilon>0$ there exists an integer $N$ so that 
	\[
		\{z\in\mathbb C:d(z,[A_f, B_f])>\epsilon\}\subset\mathbb 		C\backslash [A_{f_n}, B_{f_n}]
	\]
	and 
	\[
		|f(z)-f_n(z)|<\epsilon,\quad z\in\{z\in\mathbb 		C:d(z,[A_f,B_f])>\epsilon\},
	\]
	whenever $n\ge N$.
\end{remark}

\begin{remark}
	Note that $f=f_{\mu}$ extends as  a univalent function to 	$\mathbb C$ if and only if $\mu$ is a point mass, i.e. 	$\mu=\delta_a$ for some $a\in\mathbb R$,
	 and then $f(z)=z+a$, $z\in\mathbb C$.
\end{remark}

Denote $\Sigma$ the space of univalent functions $f:\mathbb H\to \mathbb H$ such that $f$ is the reciprocal Cauchy transform of some $\mu\in{\bf M}^U$ and endow $\Sigma$ with the metric $\rho'$ induced from $\rho$, i.e. if $f_1,f_2\in\Sigma$ and $f_1=f_{\mu_1}$, $f_2=f_{\mu_2}$, then $\rho'(f_1,f_2)=\rho(\mu_1,\mu_2)$. Let $\mathfrak P(\Sigma)$ denote the space $C([0,\infty);\Sigma)$ of continuous paths $\Psi:[0,\infty)\to\Sigma$ with the topology of uniform convergence on compact intervals induced for example by the metric
\[
	D(\Phi,\Psi)\equiv\sum_{n=1}^{\infty}\frac{1}{2^n}\cdot
	\frac{\sup_{t\in[0,n]}\rho'(\Phi(t),\Psi(t))}{1+
	\sup_{t\in[0,n]}\rho'(\Phi(t),\Psi(t))}.
\]
Then $\mathfrak P(\Sigma)$ is a polish space. If $g(t,\psi;\cdot)$ are the values of the solutions of the L\"owner equation \eqref{E:loewner} for  fixed $(t,\psi)\in[0,\infty)\times\mathfrak P(\mathbb R)$,  and where $z$ ranges over $\mathbb H\backslash{\cal K}_t$, then it follows from Proposition \ref{P:lawler} and \cite[Lemma 2]{bauer:2002} that $f\equiv g^{-1}\in\Sigma$. Finally, let $L$ be the map defined by 
\[
	\psi\in\mathfrak P(\mathbb R)\mapsto\{f(t,\psi;\cdot):\mathbb 	H\to\mathbb H, t\in[0,\infty)\}\in\mathfrak P(\Sigma).
\]
Then the probability measure $S_{\kappa}\equiv L_*{\cal W}_{\kappa}$ on $\mathfrak P(\Sigma)$ is the distribution of a stochastic L\"owner evolution with parameter $\kappa$ ($\text{SLE}_{\kappa}$).  

\begin{proposition}\label{P:weak}
	The sequence $\{L\circ S_n\}_{n=1}^{\infty}$ converges in 	distribution to $\text{SLE}_{\kappa}$, i.e. 
	\[
		L_*(S_n)_*P\Longrightarrow L_*{\cal W}_{\kappa}.
	\]
\end{proposition}

\begin{proof}
	Since $(S_n)_*P\Longrightarrow{\cal W}_{\kappa}$ as 	$n\to\infty$, it is 	enough to show that $L:\mathfrak 	P(\mathbb R)\to\mathfrak P(\Sigma)$ is continuous.
	For $t\in[0,\infty)$, $\psi\in\mathfrak P(\mathbb R)$ and 	$z\in\mathbb H$, consider the initial value problem
	\begin{equation}\label{E:initialvalue}
		\frac{\partial}{\partial s}h(s,\psi;z)
		=-\frac{2}{h(s,\psi;z)-\psi(t-s)},\ 0<s< t, \text{ and }
		h(0,\psi;z)=z.
	\end{equation}
	Then $|(\partial/\partial s) h(s,\psi;z)|\le2/\Im(h(s,\psi;z))$ and 
	\begin{equation}\label{E:imh}
		\frac{\partial}{\partial s}\Im(h(s,\psi;z))
		=\frac{2\Im(h(s,\psi;z))}{|h(s,\psi;z)-\psi(t-s)|^2}>0.
	\end{equation}
	In particular, the initial value problem has a solution for $0\le s	\le t$. Given $\psi,\varphi\in\mathfrak P(\mathbb R)$, let 	$u_1=h(\cdot,\psi;z)$, $u_2=h(\cdot,\varphi;z)$, and set
	\[
		v(s,\psi;z)=-\frac{2}{z-\psi(t-s)},\quad s\in[0,t], 		\psi\in\mathfrak P(\mathbb R),z\in\mathbb H.
	\]
	 Then
	\begin{align}
		\dot{u}_2-v(s,\psi;u_2) &=-\frac{2}{u_2-\varphi(t-s)}+
		\frac{2}{u_2-\psi(t-s)}\notag\\				&=\frac{2(\psi(t-s)-\varphi(t-s))}{(u_2-\varphi(t-s))
		(u_2-\psi(t-s))},\notag
	\end{align}
	and it follows from \eqref{E:imh} that for $s\in[0,t]$
	\begin{equation}\label{E:approxsol}
		|\dot{u}_2-v(s,\psi;u_2)|
		\le\frac{2\sup_{s\in[0,t]}|\varphi(s)-\psi(s)|}{\Im(z)^2}.
	\end{equation}
	Note also that 
	\begin{equation}\label{E:approxsol2}
		\left|\frac{\partial}{\partial z}v(s,\psi;z)\right|
		=\left|\frac{2}{(z-\psi(t-s))^2}\right|
		\le\frac{2}{\Im(z)^2}.
	\end{equation}
	Thus, by \cite[10.5.1.1]{dieudonne:1985}, if $n\in\mathbb 	Z^+$, 	then
	\[
		\sup_{s\in[0,t]}\sup_{z\in \mathbb H_{1/n}}|u_1-u_2|\le
		\left(\exp\left[2n^2t\right]-1\right)
		\sup_{s\in[0,t]}|\varphi(s)-\psi(s)|.
	\]	
	Since the initial value problem \eqref{E:initialvalue} describes 	the reverse flow to the L\"owner equation \eqref{E:loewner}, 	we 	have $f(t,\psi;\cdot)\equiv h(t,\psi;\cdot)$. Thus, for 	$n\in\mathbb Z^+$,
	\[
		\sup_{t\in[0,n]}\sup_{z\in\mathbb H_{1/n}}|f(t,\psi;z)
		-f(t,\varphi;z)|\le\left(\exp\left[2n^3 \right]-1\right)
		\sup_{t\in[0,n]}|\varphi(t)-\psi(t)|.
	\]
	Consider now the initial value problem \eqref{E:initialvalue} 	with $z=x\in\mathbb R$ and let 
	\[
		A(t,\psi)=\mathbb R\backslash\{x\in\mathbb 			R:\min_{s\in[0,t]}|h(s,\psi;x)-\psi(t-s)|>0\}.
	\]
	By continuity, $A(t,\psi)$ is connected. In fact, 	$A(t,\psi)=[A_{f(t,\psi;\cdot)}, B_{f(t,\psi;\cdot)}]$. Indeed, 	it is clear that 
	\[
		A(t,\psi)\supseteq [A_{f(t,\psi;\cdot)}, B_{f(t,\psi;\cdot)}],
	\] 
	and also
	\[
		f(t,\psi;A(t,\psi)\backslash 	[A_{f(t,\psi;\cdot)}, 		B_{f(t,\psi;\cdot)}])\subset\mathbb R.
	\]
  	Since ${\cal K}(t,\psi)=\overline{\mathbb H}\backslash 	f(t,\psi;\overline{\mathbb H}\backslash A(t,\psi))$ and 	$\overline{\mathbb H\cap {\cal K}(t,\psi)}={\cal K}(t,\psi)$  we 	have $f(t,\psi;A(t,\psi)\backslash [A_{f(t,\psi;\cdot)}, 	B_{f(t,\psi;\cdot)}])=\emptyset$.
	It now follows from Lemma \ref{L:continuity} that we can 	make the  Hausdorff distance between $A(t,\psi)$ and 	$A(t,\varphi)$ as small as we like by choosing $\psi$ close to 	$\varphi$.
\end{proof}

\begin{lemma}\label{L:continuity}
	The  Hausdorff distance between $A(t,\psi)$ and $A(t,\varphi)$ 	 is less or equal $\delta>0$ whenever	 
	\begin{equation}\label{E:bound}
		\sup_{s\in[0,t]}|\psi(s)-\varphi(s)|<\frac{\delta}{3}\wedge
		\frac{2}{3}\cdot
		\frac{\delta}{\exp(9t/(2\delta^2))-1}.
	\end{equation}
\end{lemma}

\begin{proof}
	Let $\delta>0$ be given. For $x,y\notin A(t,\psi)$ we have, 
	\begin{align}\label{E:increase}
		\frac{\partial}{\partial s}&\left[h(s,\psi;x)-\psi(t-		s)-(h(s,\psi;y)-\psi(t-s))\right]\notag\\
		&=2\frac{h(s,\psi;x)-\psi(t-s)-(h(s,\psi;y)-\psi(t-		s))}{(h(s,\psi;x)-\psi(t-	s))(h(s,\psi;y)-\psi(t-s))}
	\end{align}
	It follows that if for example $x>y>\max_{s\in[0,t]}\psi(t-	s)$  and $d(y,A(t,\psi))>0$, then $\min_{s\in[0,t]} 	|h(s,\psi;x)-\psi(t-s)|\ge x-y$. In particular, 	$d(y,A(t,\psi))\ge\delta$ implies $\min_{s\in[0,t]} 	|h(s,\psi;x)-\psi(t-s)|\ge\delta$. We will show that the latter 	together with 	\eqref{E:bound} implies that $x\notin 	A(t,\varphi)$. Then 	
	\[
		\sup_{x\in A(t,\varphi)}d(x, A(t,\psi))\le\delta.
	\]
	 By symmetry  we then also have $\sup_{x\in A(t,\psi)}d(x, 	A(t,\varphi))\le\delta$.

	If \eqref{E:bound} holds, then $|\psi(0)-\varphi(0)|<\delta/3$. 	Since also $|h(0,\psi;x)-\psi(t)|\ge\delta$, there exists 	$t_0\in(0,t]$ such that 	$\min_{s\in[0,t_0]}|h(s,\varphi;x)-\varphi(t-s)|>0$. We claim 	that we may choose $t_0=t$. For if not, then there exists 	$t_0<t'\le t$ such that $\lim_{s\nearrow 	t'}h(s,\varphi;x)=\varphi(t-t')$. Let $u=h(s,\psi;x)$ and set	$v(s,\varphi;u)=-\frac{2}{u-\varphi(t-s)}$. Then 
	\[
		\left|\frac{\partial}{\partial u}v(s,\varphi;u)\right|
		\le\frac{9}{2\delta^2},\quad s\in[0,t]
	\]
	since $|h(s,\psi;x)-\varphi(t-s)|\ge2\delta/3$. Furthermore,
	\[
		\dot{u}-v(s,\varphi;u)=\frac{2(\varphi(t-s)-\psi(t-s))}{
		 (h(s,\psi;x)
		-\psi(t-s))(h(s,\psi;x)-\varphi(t-s))}
	\]
	and so
	\[
		|\dot{u}-v(s,\varphi;u)|\le\frac{3}{\epsilon^2}\sup_{s\in[0,t]}
		|\varphi(s)-\psi(s)|.
	\]
	Again by \cite[10.5.1.1]{dieudonne:1985}
	\begin{align}\label{E:star}
		\sup_{s\in[0,t]}&|h(s,\psi;x)-h(s,\varphi;x)|\notag\\
		&\le\frac{2}{3}\sup_{s\in[0,t]}|\varphi(s)-\psi(s)|\left(\exp\left(
		\frac{9t}{2\epsilon^2}\right)-1\right)\notag\\
		&\le\frac{4}{9}\delta.
	\end{align}
	Thus, if $t_0\le t$, then $|h(t_0,\psi;x)-\psi(t-t_0)|\le 	\frac{7}{9}\delta$, a contradiction. Hence $x\notin 	A(t,\varphi)$. 	
\end{proof}
	
\begin{remark}
	Since
	\begin{align}
		\delta&\le|h(s,\psi;x)-\psi(t-s)|\notag\\
		&\le|h(s,\psi;x)-h(s,\varphi;x)|+|h(s,\varphi;x)-\varphi(t-s)|
		+|\varphi(t-s)-\psi(t-s)|\notag
	\end{align}
	the above proof together with \eqref{E:bound} and 	\eqref{E:star} implies that
	\[
		\min_{s\in[0,t]}|h(s,\varphi;x)-\varphi(t-s)|\ge\frac{2}{9}\delta.
	\]
\end{remark}

For $t\in[0,\infty)$, $\psi\in\mathfrak P(\mathbb R)$, $z\in\mathbb H$ and $n\in\mathbb Z^+$, consider the initial value problem $h_n(0,\psi;z)=z$ and 
\[
	\frac{\partial}{\partial s}h_n(s,\psi;z)=
	-\frac{2}{h_n(s,\psi;z)-\psi\left(\frac{m}{n}\right)},
\]
if $0<s<t$ and $t-s\in\left[\frac{m}{n},\frac{m+1}{n}\right)$ for some $m\in\mathbb N$. Then $|(\partial/\partial s) h_n(s,\psi;z)|\le2/\Im(h_n(s,\psi;z))$ and 
\begin{equation}\label{E:imhn}
		\frac{\partial}{\partial s}\Im(h_n(s,\psi;z))>0.
\end{equation}
In particular, the initial value problem has a solution for $0\le s\le t$. Proposition \ref{P:lawler} extends to piecewise continuous $\psi$ and thus $f_n(t,\psi;\cdot)\equiv h_n(t,\psi;\cdot)\in\Sigma$. Let $L_n$ be the map defined by
\[
	\psi\in\mathfrak P(\mathbb R)\mapsto\{f_n(t,\psi;\cdot):
	\mathbb H\to\mathbb H, t\in[0,\infty)\}\in\mathfrak P(\Sigma).
\]
We can consider the family of random variables $\{(L_n\circ S_n)(\frac{m}{n})\}_{m=0}^{\infty}$ as a random walk on $\Sigma$ as follows. 
For $a\in\mathbb R$ let $r_n(a;\cdot)$ be the conformal map given by
\[
	z\in\mathbb H\mapsto r_n(a;z)
	=a+\sqrt{(z-a)^2-\frac{4}{n}}\in\mathbb H.
\]
Then
\begin{equation}\label{E:rangefan}
	r_n(a;\mathbb H)=\mathbb H\backslash\{z\in\mathbb 	H:\Re(z)=a\text{ and }\Im(z)\in[0,\frac{2}{\sqrt{n}}]\}.
\end{equation}
For $n\in\mathbb Z^+$, $\omega\in\Omega$, and $z\in\mathbb H$ set $D_n(0,\omega;z)=z$ and define inductively
\[
	D_n(m,\omega;z)= 	D_n\left(m-1,\omega;r_n\left(S_n\left(\frac{m-1}{n},\omega 	\right);z\right)\right),
\]
if $m>0$. Then, for every $n\in\mathbb Z^+$ and $\omega\in\Omega$, $\{D_n(m,\omega;\cdot)\}_{m=0}^{\infty}$ is a family of conformal maps from $\mathbb H$ into $\mathbb H$ and, 
\begin{equation}\label{E:shrinkingimage}
	D_n(m,\omega;\mathbb H)\supsetneq D_n(m+1,\omega;
	\mathbb H),\quad \text{ for every } m\in\mathbb N.
\end{equation}
In fact, $r_n(a;z)$ is the solution at time $t=1/n$ of the initial value problem  $(\partial/\partial s)h(s,a;z)=-2/(h(s,a;z)-a)$, $h(0,a;z)=z$. Thus $r_n(a;\cdot)\in\Sigma$ for every $a\in\mathbb R$,  $D_n(m,\omega;\cdot)\in\Sigma$ for every $m\in\mathbb  N$ and $\omega\in\Omega$, and finally
\[
	D_n(m,\omega;z)
	=f_n\left(\frac{m}{n},S_n(\cdot,\omega);z\right).
\]
By boundary correspondence, $D_n(m,\omega;\cdot)$ maps the real axis to a finite number of Jordan arcs. All prime ends are of the first kind and hence $D_n(m,\omega;\cdot)$ extends continuously to $\bar{\mathbb H}$, see \cite[Theorem 2.21]{markushevich:1967}. 

\begin{theorem}\label{T:dle}
	  The sequence $\{L_n\circ S_n\}_{n=1}^{\infty}$ converges in 	distribution to $\text{SLE}_{\kappa}$, i.e. 
	\[
		(L_n)_*(S_n)_*P\Longrightarrow L_*{\cal W}_{\kappa}.
	\]
\end{theorem}

\begin{proof}
	With the notation from above we have
	\begin{align}
		\dot{h}_n-v(s,\psi;h_n) &=
		-\frac{2}{h_n-\psi\left(\frac{m}{n}\right)}+
		\frac{2}{h_n-\psi(t-s)}\notag\\				&=\frac{2\left(\psi(t-s)
		-\psi\left(\frac{m}{n}\right)\right)}{\left(h_n
		-\psi\left(\frac{m}{n}\right)\right)(h_n-\psi(t-s))},\notag
	\end{align}
	and it follows from \eqref{E:imhn} that for $s\in[0,t]$
	\begin{equation}\label{E:approxsol}
		|\dot{h}_n-v(s,\psi;h_n)|
		\le\frac{2\rho(n,t;\psi)}{\Im(z)^2},
	\end{equation}
	where $\rho(n,t;\psi)\equiv\sup\{|\psi(r)-\psi(s)|:0\le s<r\le t
	\text{ with }r-s\le\frac{1}{n}\}$ is the modulus of continuity of 	$\psi$. Thus, from \cite[10.5.1.1]{dieudonne:1985}, if 	$N\in\mathbb Z^+$,
	\[
		\sup_{s\in[0,t]}\sup_{z\in\mathbb H_{1/N}} 		|h_n(s,\psi;z)-h(s,\psi;z)|\le\left(\exp[2N^2t]-1\right)
		\rho(n,t;\psi).
	\]
	Similarly, the proof of Lemma \ref{L:continuity} extends to 	show that that $A(t,\psi_n)\to A(t,\psi)$ in the Hausdorff 	distance, as $n\to\infty$, where $\psi_n$ is defined by 	$\psi_n(t-s)=\psi(m/n)$ if $t-s\in[m/n,(m+1)/n)$, and where we 	now define $A(t,\psi)$ as the convex closure of $\mathbb 	R\backslash\{x\in\mathbb 			R:\min_{s\in[0,t]}|h(s,\psi;x)-\psi(t-s)|>0\}$. 
	It follows that, for each $\psi\in\mathfrak P(\mathbb R)$, 	$D(f_n(\cdot,\psi;\cdot),f(\cdot,\psi;\cdot))\to0$ as $n\to\infty$. 	In particular, $D(L_n\circ S_n,L\circ S_n)\to0$ in probability as 	$n\to\infty$ and by the principle of accompanying laws, 
	\cite[3.1.14]{stroock:1993}, and Proposition \ref{P:weak} we 	get
	\[
 		(L_n)_*(S_n)_*P\Longrightarrow L_*{\cal W}_{\kappa}.
	\]
\end{proof}
\section{Properties of  discrete L\"owner evolution} \label{S:dle}

For all $n\in\mathbb Z^+$,
\[
	r_n(a;z)=\frac{1}{\sqrt{n}}r_1\left(\sqrt{n}a;\sqrt{n}z\right),
	\quad a\in\mathbb R, z\in\mathbb H,
\]
and  
\[
	D_n\left(0,\omega;\frac{z}{\sqrt{n}}\right)=
	\frac{1}{\sqrt{n}}D_1(0,\omega;z),\quad \omega\in\Omega,
	z\in\mathbb H.
\]
Using 
\[
	\sqrt{n}S_n\left(\frac{m}{n},\omega\right)=\sum_{k=1}^m 	X_k(\omega)=S_1(m,\omega),\quad\omega\in\Omega,
	m\in\mathbb N,
\]
it follows by induction that
\[
	D_n\left(m,\omega;\frac{z}{\sqrt{n}}\right)=
	\frac{1}{\sqrt{n}}D_1(m,\omega;z),\quad \omega\in\Omega,
	z\in\mathbb H,m\in\mathbb N.
\]
Thus to study the families $\{D_n(m)\}_{m=0}^{\infty}$ we may as well restrict to $n=1$. 

Writing $D$, $S$, $r$ for $D_1$, $S_1$, and $r_1$, respectively, $\{D(m)\}_{m=0}^{\infty}$ is defined by
\begin{equation}\label{E:d0}
    D(0,\omega;z)=z,\quad\omega\in\Omega,z\in\mathbb H,
\end{equation}
and 
\begin{equation}\label{E:dm}
	D(m,\omega;z)=D(m-1,\omega;r[S(m-1,\omega);z]),\quad
	\omega\in\Omega,z\in\mathbb H,
\end{equation}
if $m\in \mathbb Z^+$ .  

\begin{theorem}\label{T:iid}
	For $(m,n)\in\mathbb N^2$ such that $m\le n$, and 	$\omega\in\Omega$, define conformal maps 	$H(m,n,\omega;\cdot):\mathbb H\to\mathbb H$ by
	\[
		H(m,n,\omega;\cdot)=D(m,\omega;\cdot)^{-1}\circ
		D(n,\omega;\cdot)
	\]
	and set
	\[
		\tilde{H}(m,n,\omega;z)
		=H(m,n,\omega;z+S(m,\omega))-S(m,\omega).
	\]
	Then the family $\{\tilde{H}(m,n)\}_{n=m}^{\infty}$ is 	independent of the family $\{D(k)\}_{k=0}^{m+1}$. 
	Furthermore, if the random variables $\{X_k\}_{k=0}^{\infty}$  	are  identically distributed, then the distribution of the random 	variable 
	$\omega\in\Omega\mapsto\tilde{H}(m,n,\omega;\cdot)
	\in\Sigma$ under $P$ is the same as the 	distribution of $\omega\in\Omega\mapsto D(n-	m,\omega;\cdot)\in\Sigma$ under $P$.  Finally, if  $X_k$ is 	symmetric for each 	$k\in\mathbb N$, then, for each 	$m\in\mathbb N$, $\omega\in\Omega\mapsto 	D(m,\omega;\cdot)\in\Sigma$ 	and 	$\omega\in\Omega\mapsto \chi\circ 	D(m,\omega;\cdot)\circ\chi\in\Sigma$ have the same 	distribution under $P$, where
	the map $\chi$ is given by
	\[
		x+\sqrt{-1}y\in\mathbb H\mapsto
		-x+\sqrt{-1}y\in\mathbb 		H.
	\] 
\end{theorem}

\begin{proof}
	For the first statement, note that the case $n=m$ is trivial and 	consider the case $n>m$. From \eqref{E:dm} and induction on 	$n$ it follows that
	\begin{equation}\label{E:expressdm}
			D(n;z)=D\left(m;r\left[S(m);r\left[S(m+1);\dots; 			r[S(n-1);z]\cdots\right]\right]\right).
	\end{equation}
	In particular, for any $m\in\mathbb Z^+$, $D(m)$ is 	$\sigma(X_1,\dots,X_{m-1})$-measurable. Furthermore, we 	now get
	\[
		\tilde{H}(m,n;z)=r\left[S(m);r\left[S(m+1);\dots; 			r[S(n-1);z+S(m)]\cdots\right]\right]-S(m).
	\]
	Applying repeatedly the identity $r(a-b;z-b)+b=r(a;z)$ gives
	\begin{align}\label{E:increment}
		\tilde{H}&(m,n;z)\notag\\
		&=r[S(m)-S(m);r[S(m+1)-S(m);\dots;r[S(n-1)
		-S(m);z]\cdots]].
	\end{align}
	This implies the first statement because the random variables 	$\{X_k\}_{k=0}^{\infty}$ are mutually independent. The 	expression \eqref{E:increment} also implies the second 	statement, under the assumption that the $X_k$'s are 	identically distributed.
	Regarding the third statement, note that $r(a;\cdot)\circ\chi=
	\chi\circ r(-a;\cdot)$ and that $\chi^{-1}=\chi$. Thus it follows 	from \eqref{E:expressdm} that
	\[
		\chi\circ D(m,\omega;\cdot)\circ\chi
		=\tilde{D}(m,\omega;\cdot),
	\]
	where $\tilde{D}(m)$ is defined as $D(m)$ in \eqref{E:dm}, 	with $-S(m)$ in place of $S(m)$, $m\in\mathbb N$. The 	symmetry of the $X_k$'s implies the symmetry of the $S(m)$'s 	and so the distributions of $\omega\in\Omega\mapsto 	D(m,\omega)\in\Sigma$ and 	$\omega\in\Omega\mapsto\tilde{D}(m,\omega)\in\Sigma$ 	under $P$ are equal. 
\end{proof}

\begin{remark} The above theorem gives the discrete version of 	corresponding results for SLE, see \cite{lawler:2001}. The first 	statement shows that $\{D(m)\}_{m=0}^{\infty}$ has, 	up to a shift, independent increments relative to 	composition of maps, the second statement shows that the  	shifted increments are stationary (assuming the $X_k$'s are 	identically distributed) and the third statement is a kind of 	reflection symmetry of $D(m)$ (for symmetric $X_k$'s) . Note 	that if we look at the images of the maps 	$D(m,\omega;\cdot)$, then we find that 	$\omega\mapsto D(m,\omega;\mathbb H)$ and 	$\omega\mapsto \chi(D(m,\omega;\mathbb H)$ have the same 	distribution under $P$ on a suitably defined space of domains 	in the upper half-plane, because $\chi(\mathbb H)=\mathbb 	H$. Equivalently, the distributions of the ``hulls'' $\mathbb 	H\backslash  D(m;\mathbb H)$ is invariant under 	$\chi$. This is the reflection symmetry statement in 	\cite{lawler:2001}. In fact, the weak convergence of the 	increments and the continuity of the map $L$ imply that the 	continuous results for SLE can be deduced directly from 	Theorem \ref{T:iid}.
\end{remark}

Assume now that $\{X_n\}_{n=0}^{\infty}$ is a sequence of independent and identically distributed random variables of mean-value $0$ and variance $\kappa\ge0$ on a probability space $(\Omega,{\cal F},P)$.
In that case we call the family $\{D(m)\}_{m=0}^{\infty}$ a {\it discrete L\"owner evolution with parameter} $\kappa$ ($\text{DLE}_{\kappa}$). 

When the $X_n'$s are Bernoulli random variables, i.e.
\[
	P(X_n=\sqrt{\kappa})=P(X_n=-\sqrt{\kappa})=\frac{1}2,\quad
	n\in\mathbb N,
\]  
then the corresponding discrete L\"owner evolution has a trivial ``phase transition'' at $\kappa=4$.

\begin{proposition}\label{P:connect}
	Let $\{D(m)\}_{m=0}^{\infty}$  be a discrete L\"owner 	evolution with parameter $\kappa$ driven by a sequence of 	Bernoulli random variables as above. If $\kappa\le4$, then 	$\overline{\mathbb H\backslash D(m,\omega;\mathbb H)}$ is 	connected in $\overline{\mathbb H}$ for all 	$\omega\in\Omega$ and $m\in \mathbb N$.
	If $\kappa>4$, then $\overline{\mathbb H\backslash 	D(m,\omega;\mathbb H)}$ is not connected in 	$\overline{\mathbb H}$, for all $\omega\in\Omega$ and 	$m\ge2$.
\end{proposition} 

\begin{proof}
	This follows immediately by considering the composition of 	maps	$r(0;\cdot)\circ r(\sqrt{\kappa};\cdot)$. The closure of 	the complement of the image of $\mathbb H$ under this map 	is connected in $\overline{\mathbb H}$ if and only if 	$\kappa\le4$.  
\end{proof}

Graphically, for $\kappa\le4$ the omitted set, i.e. $\overline{\mathbb H\backslash D(m,\omega;\mathbb H)}$, is a single tree made up of $m$ curvy branches. The tree grows one branch at each step. Orient the branches in the direction of the root and label the end-point closest to the root ``bottom'' and the other end-point ``top''. If $0<\kappa<4$, then the $(m+1)$st branch ``branches off'' the $m$th branch somewhere between the $m$th branch's top and bottom. We call the segment of the $m$th branch between the branch-point to the $(m+1)$st branch and the top of  of the $m$th branch the overshoot. For $\kappa=4$ the branch-point is at the bottom and the tree looks like a bushel, all branches emanating from the point $z=0$, while for $\kappa=0$ the branch point is at the top and the tree degenerates to a vertical line segment in the closed upper half-plane beginning at $z=0$. As $\kappa$ decreases from $4$ to $0$ the branch-point increases from bottom to top. Using the orientation towards the root, the $(m+1)$st branch branches off to the right of the $m$th branch if $X_m>0$, and to the left if $X_m<0$.

If $\kappa>4$, $\overline{\mathbb H\backslash D(m,\omega;\mathbb H)}$ consists of $m$ branches forming at least $\min(m,2)$ trees and at most $m$ trees. The latter will be the case for instance if the driving random walk makes all of its first $m-1$ steps in one direction, while the former picture emerges if the walk changes direction at every step. Typically, for large $m$ neither will be the case and it would be interesting for example to calculate the expected number of trees, or the distribution of the distance of the roots of neighboring trees. For example, by first letting the random walk alternate directions for a long time and then stepping only in one direction for a long time, it is easy to see that roots may be spaced arbitrarily far apart.

If $X_n$ is centered and of variance $\kappa$ but not necessarily a Bernoulli random  variable then the above picture should still be approximately right. Of course, even for $\kappa\le4$ we may now get several trees. But their number or spacing should be small as $m\to\infty$ compared to the case when $\kappa>4$.

The phase transition for $SLE$ at $\kappa=4$ is the fact that ${\cal K}(t)$ is a simple curve for $\kappa\le4$, $P-$a.s., and that it is not a simple curve for $\kappa>4$, $P-$a.s, \cite{rohde.schramm:2001}. Thus, in the scaling limit, the overshoots disappear, creating a simple curve if $\kappa\le4$. For $\kappa\ge4$, the disjoint trees become connected in the scaling limit (if they are too small, some might also disappear). 

We now study a question related to this phase transition following ideas in \cite{lawler:2001}.

Let $g_n(m)=(D_n(m))^{-1}$, $n\in\mathbb Z^+$, $m\in\mathbb N$. Then
\[
	g_n(m,\omega;\cdot)=\left(r_n\left[S_n\left(
	\frac{m-1}{n},\omega\right);\cdot\right]\right)^{-1}\circ 
	g_n(m-1,\omega;\cdot),
\]
that is
\[
	g_n(m,\omega;\cdot)=S_n\left(\frac{m-1}n,\omega\right)+
	\sqrt{\left[g_n(m-1,\omega;\cdot)-S_n\left(\frac{m-1}n,
	\omega\right)\right]^2+\frac{4}{n}}.
\]  
In particular, if we set 
\[
	Y_n(m,\omega;\cdot)=\frac{g_n(m,\omega;
	\cdot)-S_n\left(\frac{m-1}{n},\omega\right)}{\sqrt{\kappa}},
\]
then
\[
	Y_n(m,\omega;\cdot)=\sqrt{\left(Y_n(m-1,\omega;\cdot)
	-\frac{X_m}{\sqrt{n\kappa}}\right)^2+\frac{4}{n\kappa}}.
\]
Note that $Y_1(m,\omega;\cdot)=\sqrt{n} Y_n(m,\omega;\cdot)$ and $X_m'\equiv X_m/\sqrt{\kappa}$ is centered with variance 1. For $z=x\in\mathbb R\backslash\{0\}$ set  $Y_0=x$ and $Y_m=Y_1(m,\omega;x)$. Then $\{Y_m\}_{m=0}^{\infty}$ is a Markov chain satisfying the evolution equation
\[
	Y_m=\sqrt{(Y_{m-1}-X_m')^2+\frac{4}{\kappa}}, \quad 	m\in\mathbb Z^+,
\]
or equivalently
\[
	Y_m^2-Y_{m-1}^2=-2Y_{m-1} X_m'+(X_m')^2
	+\frac{4}{\kappa},\quad m\in\mathbb Z^+.
\]

\begin{theorem}\label{T:Bessel}
	The Markov chain $\{Y_m\}_{m=0}^{\infty}$
	is transient if $\kappa<4$, and it is recurrent if $\kappa>4$. If 	the moment generating function of $X_1'$ has a positive 	radius of convergence, then the chain $\{Y_m\}$ is recurrent 	for $\kappa=4$.    
\end{theorem}

\begin{proof}
	It is easy to see that $\limsup_{m\to\infty}Y_m=+\infty$, 	$P$-a.s.
	Using Taylor series and a cutoff for $X_m'$ if necessary we 	see that
	\begin{align}
		\lim_{y\to\infty}&2y\ \mathbb E^P[Y_m-		Y_{m-1}|Y_{m-1}=y]\notag\\
		&=\lim_{y\to\infty}2y^2\ \mathbb 			E^P\left[\sqrt{\left(1-\frac{X_m'}{y}\right)^2
		+\frac{4}{y^2\kappa}}-1\right]  \notag\\
		&=\lim_{y\to\infty}2y^2\ \mathbb E^P\left[-\frac{X_m'}y
		+\frac{(X_m')^2}{2y^2}+\frac{2}{y^2\kappa}-
		\frac{(X_m')^2}{2y^2}\right]=\frac{4}{\kappa}.\notag
	\end{align}
	Under additional moment assumptions, the convergence is at 	least of order $O(1/y)$. Furthermore,
	\begin{align}
		\mathbb E^P[&(Y_m- Y_{m-1})^2|Y_{m-1}=y]\notag\\
		&=
		\mathbb E^P[-2Y_{m-1}X_m'+(X_m')^2+4/\kappa
		-2Y_{m-1}(Y_m-Y_{m-1})|Y_{m-1}=y]\notag\\
		&=1+\frac{4}{\kappa}-2y\ \mathbb E^P
		[Y_m-Y_{m-1}|Y_{m-1}=y]\to1,\notag
	\end{align}
	as $y\to\infty$. Thus, by \cite[Theorem 3.2]{lamperti:1960}, the 	result follows.
\end{proof}

\appendix
\section{Monotonic Independence and L\"owner Map}

The following definition is taken from \cite{muraki:2001}. Let $({\cal A},\phi)$ be a $C^*$-probability space consisting of a unital $C^*$-algebra $\cal A$ and a state $\phi$ over $\cal A$. The elements $X$ of $\cal A$ are called random variables and $\phi(X)$ their expectation.

\begin{definition}
	A family $\{X_i\}_{i\in I}\subset\cal A$ of random variables on 	$({\cal A},\phi)$ with totally ordered index set $I$ is said to be 
	monotonically independent with respect to a state $\phi$ if the 	following two conditions are satisfied.  
	
	$(a)$ Whenever $i<j$, $k<j$, and $p\in\mathbb N$, then
	\[
		X_i X_j^p X_k=\phi(X_j^p) X_i X_k.
	\]
		
	$(b)$ Whenever $i_m>\cdots>i_1>i$, $j_n>\cdots>j_1>i$, and 
	$p,p_k,q_l\in\mathbb N$, 
	then
	\begin{align}
		\phi(X_{i_m}^{p_m}&\cdots X_{i_1}^{p_1}X_i^p 		X_{j_1}^{q_1}\cdots X_{j_n}^{q_n})\notag\\
		&=\phi(X_{i_m}^{p_m})\cdots\phi(X_{i_1}^{p_1})\phi(X_i^p)
		\phi(X_{j_1}^{q_1})\cdots\phi(X_{j_n}^{q_n}).\notag
	\end{align}
\end{definition}

\begin{theorem}\cite{muraki:2000}
	Let $X_1, X_2,\dots, X_n\in\cal A$ be monotonically 	independent self-adjoint random variables on $({\cal A},\phi)$, 	in the natural order of  $\{1,2,\dots,n\}$. If $f_{X_k}:\mathbb 	H\to\mathbb H$ denotes the reciprocal Cauchy transform of 	the distribution of $X_k$, for $1\le k\le n$, then
	\[
		f_{X_1+X_2+\cdots+X_n}=f_{X_1}\circ f_{X_2}\circ\cdots
		\circ f_{X_n}.
	\]
\end{theorem}   	

Define for a pair of probability measures $\mu$, $\nu$ on $\mathbb R$ the monotonic convolution $\lambda$ of $\mu$ and $\nu$,  denoted by $\lambda=\mu\rhd\nu$, as the unique probability measure $\lambda$ satisfying $f_{\lambda}(z)=f_{\mu}(f_{\nu}(z))$. T

\begin{corollary}
	For $\psi\in\mathfrak P(\mathbb R)$ let $f(t,\psi)=g^{-1}(t,\psi)$ 	be the solution to the L\"owner equation \eqref{E:loewner}. For 	$0\le s\le t$ set $f_{s,t}=g_s\circ f_t$. Then 	$f_{s,t}=f_{\mu_{s,t}}$ for a unique probability measure 	$\mu_{s,t}$, and, for $r\le s\le t$,
	\[
		\mu_{r,s}\rhd\mu_{s,t}=\mu_{r,t}.
	\]
	Similarly, for $\omega\in\Omega$ and $0\le m\le n$, 	$H(m,n,\omega)$ is the reciprocal Cauchy transform of a 	unique probability measure $\mu_{m,n}$ , and if $l\le m\le n$,
	then
	\[
		\mu_{l,m}\rhd\mu_{m,n}=\mu_{l,n}.
	\]
\end{corollary}

Thus $\{f(t,\psi):t\in[0,\infty)\}$ and $\{D(m,\omega):m\in\mathbb N\}$ correspond to  monotonically independent increment processes in some noncommutative probability space $({\cal A},\phi)$.  In fact, the ``building blocks'' for our discrete L\"owner evolution, the functions $r_n(a;z)=a+\sqrt{(z-a)^2-4/n}$, are the reciprocal Cauchy transforms of some well known distributions: the arcsine distribution supported in $(-2/\sqrt{n},2/\sqrt{n})$ if $a=0$, and a deformation of the arcsine distribution if $a\neq0$. Note that the arcsine distribution plays for monotonic convolution the role the Gaussian distribution plays for ``classical convolution.'' For example, the monotonic central limit theorem establishes convergence to an arcsine distribution.  




\end{document}